\documentclass[12pt,twoside]{article}
\usepackage{amsfonts}
\usepackage{amssymb,amsmath}
\begin{document}
\begin{center} \textbf{Mathematical Structures Defined by Identities II} \vspace{12pt}
\\Constantin M. Petridi
\\ cpetridi@hotmail.com
\\
\vspace{0.5cm}
\end{center}
\leftskip = 1cm \rightskip = 1cm\footnotesize{\textbf{Abstract} In
our paper arXiv: math.RA/0110333 v1 Oct 2001 we showed that the
number of algebras defined by a binary operation satisfying a
formally irreducible identity between two $n$-iterates is
$O\bigl( e^{-n/16}S_{n}^{2}\bigr)$ for $n\rightarrow\infty$,
$S_{n}$ being the $n$th-Catalan number. This was proved by using
exclusively the series of tableaux $A_{n}$. By using also the
series of tableaux $B_{n}$, we now sharpen this result to
$O\Bigl(|\frac{n+2}{n}\;e^{-n/16}-\frac{2}{n}|S_{n}^{2}\Bigr).$}
\par
\leftskip = 0cm \rightskip = 0cm \normalsize{The exposition
follows, in abbreviated form, the outline of above arXiv paper,
denoted by MS, to which we refer for explanation of concepts and
symbols.\\
\par
\textbf{1.} Since tableau $A_{n}$ has n-lines and tableau $B_{n}$
has 2 lines, tableau $A_{n}\oplus B_{n}$ has $n+2$ lines. The
relations of MS 2.2 regarding the number of their (lines) common
elements, have to be unchanged as follows, so that we again have
for $k=1,2,\dots,n+2$
$$|L_{i_{1}}\cap L_{i_{2}}\cap \dots \cap L_{i_{k}}| =
\left\{\begin{matrix} S_{n-1}\hfill&if \;\; i_{1}=i_{2}=\dots
=i_{k}\hfill &\cr 0 \hfill &\text{if at least one }
|i_{1}-i_{2}|,\dots,|i_{k-1}-i_{k}|,\hfill\cr &\text{all taken
mod\, n,\;\;is equal to 1}\hfill \cr S_{n-k}&
\text{otherwise.}\hfill
\end{matrix}\right. $$
For example, for $n=6,\;k=2,$ there are 8 lines
$L_{1},\;L_{2},\;\dots,L_{8}$ in tableau $A_{6}\oplus B_{6}$. The
$8\times 8$ table $(|L_{i}\cap L_{j}|)$ looks
as follows (only the entries on and above the diagonal are shown since $|L_{i}\cap L_{j}|=|L_{j}\cap L_{i}|$.) \\
\scriptsize{
$$\left. \begin{matrix} &L_{1}&L_{2}&L_{3}&L_{4}&L_{5}&L_{6}&L_{7}&L_{8}\cr
& & & & & & & & \cr L_{1}
&S_{6}&0&S_{5}&S_{5}&S_{5}&S_{5}&S_{5}&0\cr L_{2} &
&S_{6}&0&S_{5}&S_{5}&S_{5}&S_{5}&S_{5}\cr L_{3} & &
&S_{6}&0&S_{5}&S_{5}&S_{5}&S_{5}\cr L_{4} & & &
&S_{6}&0&S_{5}&S_{5}&S_{5}\cr L_{5} & & & &
&S_{6}&0&S_{5}&S_{5}\cr L_{6} & & & & & &S_{6}&0&S_{5}\cr L_{7} &
& & & & & &S_{6}&0\cr L_{8} & & & & & & & &S_{6}\cr
\end{matrix}\right. \;\;\;=\;\;\;
\left. \begin{matrix}
&L_{1}&L_{2}&L_{3}&L_{4}&L_{5}&L_{6}&L_{7}&L_{8}\cr & & & & & & &
& & \cr L_{1}&42&0&14&14&14&14&14&0\cr L_{2}&
&42&0&14&14&14&14&14\cr L_{3}& & &42&0&14&14&14&14\cr L_{4}& & &
&42&0&14&14&14\cr L_{5}& & & & &42&0&14&14\cr L_{6}& & & & &
&42&0&14\cr L_{7}& & & & & & &42&0\cr L_{8}& & & & & & & &42\cr
\end{matrix}\right.$$}
\par
\par
\normalsize \textbf{2}. The multiplicity $M(J_{i}^{n})$ of an
$n$-interate $J_{i}^{n}$ is the number of times the iterate
occurs in tableau $A_{n}\oplus B_{n}$.The number of $n$-iterates
with same multiplicity $k$, is now denoted by
$T_{n,k}^{A_{n}\oplus B_{n}}$. As stated at the end of MS  2.6
this number has been found to be, for $k\geq 1$
%---------------------------------------equation 1-------------------------------------
$$\hspace{3.4cm} T_{n,k}^{A_{n}\oplus B_{n}}=T_{n,k}+2T_{n-1,k-1}-2T_{n-1,k},\hspace{3.4cm}   (1)$$
where $T_{n,k}$ are the corresponding numbers with regard to
tableau $A_{n}$, i.e.
%--------------------------------------equation 2--------------------------------------
$$\hspace{2cm}T_{n,k}=2^{n-2k+1}\binom{n-1}{2k-2}S_{k-1},\;\;\; k=1,2,\dots,\hspace{2cm} \hspace{1.5cm}(2)$$
%----------------------------end page 1---------------------------------------------
Calculating the terms of the left side of (1), we have from (2)
$$T_{n,k}=2^{n-2k+1}\binom{n-1}{2k-2}S_{k-1}=2^{n-2k+1}
\frac{(n-1)\dots (n-2k+2)}{(2k-2)!}S_{k-1}$$
$$2T_{n-1,k-1}=2^{n-2k+3}\binom{n-2}{2k-1}S_{k-2}=2^{n-2k+3}
\frac{(n-2)\dots (n-2k+3)}{(2k-4)!}S_{k-2}$$
$$-2T_{n-1,k}=-2^{n-2k+1}\binom{n-2}{2k-2}S_{k-1}=-2^{n-2k+1}
\frac{(n-2)\dots (n-2k+1)}{(2k-2)!}S_{k-1}.$$
Adding and using
the recursion $S_{k}=2\frac{2k-1}{k+1}S_{k-1}$ for the Catalan
numbers we get
%-----------------------------------------------equation 3----------------------------------
$$T_{n,k}^{A_{n}\oplus B_{n}}=2^{n-2k+1}\frac{(n+2)\dots (n-2k+3)}{(2k-4)!}\Bigg\{
\frac{(n-1)(n-2k+1)}{(2k-3)(2k-2)}S_{k-1}+$$
$$4S_{k-2}-\frac{(n-2k+2)(n-2k+1)}{(2k-3)(2k-2)}S_{k-1}\Bigg\}$$
$$=2^{n-2k+1}\binom{n-2}{2k-4}\Bigg\{
\frac{n-2k+2}{(2k-3)(2k-2)}(n-1-n+2k-1)S_{k-1}+4S_{k-2}\Bigg\}$$
$$=2^{n-2k+1}\binom{n-2}{2k-4}\Bigg\{
\frac{n-2k+2}{2k-3}S_{k-1}+4S_{k-2}\Bigg\}$$
$$=2^{n-2k+1}\binom{n-2}{2k-4}\Bigg\{
\frac{n-2k+2}{k}+2\Bigg\}2S_{k-2}$$
$$=2^{n-2k+2}\binom{n-2}{2k-4}\frac{n+2}{k}S_{k-2}. $$
But from (2) we have that
$$T_{n-1,k-1}=2^{n-2k+2}\binom{n-2}{2k-4}S_{k-2}$$
so that finally
%-------------------------------------------equation 4---------------------------------------
$$\hspace{4cm}T_{n,k}^{A_{n}\oplus B_{n}}=\frac{n+2}{k}T_{n-1,k-1}.\hspace{4cm}(3)$$
\par
\textbf{3.} Formal reducibility of an identity
$J_{i}^{n}=J_{j}^{n}$  and incidence matrix relative to tableau
$A_{n}\oplus B_{n}$ are defined in the same way as for tableau
$A_{n}$. The number of formally \textit{reducible} identities
$J_{i}^{n}=J_{j}^{n}$ of order $n$, which we denote by
$I_{n}^{A_{n}\oplus B_{n}}$, to distinguish it from $I_{n}$
relative to tableau $A_{n}$, is given by
$$I_{n}^{A_{n}\oplus B_{n}}=\sum_{1\leq i,j\leq S_{n}}\delta
(J_{i}^{n},J_{j}^{n}),$$
%---------------------- end page 2 --------------------
where
$$\delta (J_{i}^{n},J_{j}^{n}) = \left\{\begin{matrix}1 & \text{if} &
J_{i}^{n}=J_{j}^{n} & \text{formally reducible}\hfill \cr 0 &
\text{if} & J_{i}^{n}=J_{j}^{n} & \text{formally
irreducible}\hfill
\end{matrix}\right.$$
\par
As an example, we display the incidence
matrices relative to tableaux $A_{3}$ and $A_{3}\oplus B_{3}$,
which clearly shows, as expected that $I_{3}^{A_{3}\oplus
B_{3}}>I_{3}.$
$$\hspace{2cm}A_{3} \hspace{6cm}A_{3}\oplus B_{3}\hspace{1.5cm}$$
$$\left. \begin{matrix} & J_{1}^{3}
&J_{2}^{3}&J_{3}^{3}&J_{4}^{3}&J_{5}^{3}&\sum_{i}1\cr
J_{1}^{3}&1&1&0&0&0&2 \cr J_{2}^{3}&1&1&0&0&1& 3\cr
J_{3}^{3}&0&0&1&1&0& 2\cr J_{4}^{3}&0&0&1&1&0& 2\cr
J_{5}^{3}&0&1&0&0&1& 2\cr &&&&&&-\cr &&&&I_{3}&=&11
\end{matrix}\right.\hspace{1cm}
\left. \begin{matrix}& J_{1}^{3}
&J_{2}^{3}&J_{3}^{3}&J_{4}^{3}&J_{5}^{3}&\sum_{i}1\cr
J_{1}^{3}&1&1&1&0&0&3 \cr J_{2}^{3}&1&1&0&0&1& 3\cr
J_{3}^{3}&1&0&1&1&0& 3\cr J_{4}^{3}&0&0&1&1&1& 3\cr
J_{5}^{3}&0&1&0&1&1& 3\cr &&&&&&-\cr &&&&I_{3}^{A_{3}\oplus
B_{3}}&=&15
\end{matrix}\right.$$

\par
For $n=4$ and $n=5$ the corresponding findings are
$$\left.\begin{matrix}
n=4 &&&I_{4}=88\hfill &&& I_{4}^{A_{4}\oplus B_{4}}=116\hfill\cr
n=5 &&& I_{5}=834 &&& I_{5}^{A_{4}\oplus B_{4}}=1050 \end{matrix}
\right.$$
\par
\par
\textbf{4.} The arguments which led us to establish the
fundamental Theorem of MS 2.4 can be applied verbatim, resulting
in
$$\hspace{2.5cm} \sum_{j=1}^{S_{n}}\delta
(J_{i}^{n},J_{j}^{n})=\sum_{\nu=1}^{M(J_{i}^{n})}(-1^{\nu-1})\binom{M(J_{i}^{n})}{\nu}S_{n-\nu},\hspace{2.5cm}(4)$$
where $M(J_{i}^{n}$ is the multiplicity of $J_{i}^{n}$ in tableau
$A_{n}\oplus B_{n}$ . (4) means that the number of formally
reducible identities in the line $L_{i}$ of tableau
$I_{n}^{A_{n}\oplus B_{n}}$ does  \textit{not} depend on
$J_{i}^{n}$ but only on the multiplicity of $M(J_{i}^{n})$, as was
the case in MS 2.4, when solely the series of tableaux $A_{n}$ was
used. Since $\binom{M(J_{i}^{n})}{\nu}=0$, for $\nu>M(J_{i}^{n})$
we can forget the upper limit $M(J_{i}^{n})$ for the index $\nu$
and write instead
$$\hspace{2.7cm}\sum_{j=1}^{S_{n}}\delta
(J_{i}^{n},J_{j}^{n})=\sum_{\nu=0}^{\infty}(-1^{\nu-1})\binom{M(J_{i}^{n})}{\nu}S_{n-\nu}.\hspace{2.7cm}(5)$$
This convention will be used for all finite series of the form
$\sum_{\nu=n}^{N}\binom{f(N)}{\nu}c_{\nu}$, $f(x)$ a positive
arithmetic function which from now on will be written as
$\sum_{\nu=n}^{\infty}\binom{f(N)}{\nu}c_{\nu}$, since
$\binom{f(N)}{\nu}=0$ for $\nu > f(N)$.
\par
\par
 %----------------------------------end page 3 ------------------------------
\textbf{5.} Following identity from MS 2.5 for various values of
the indices $n$ and $k$ will be needed in the sequel
$$\sum_{\nu=0}^{[ \frac{n+1}{2}]-k}\binom{k+\nu}{k}T_{n,k+\nu}=
 \binom{n-k+1}{k}S_{n-k},$$
which, because of said convention, will be written
 $$\hspace{3cm}\sum_{\nu=0}^{\infty}\binom{k+\nu}{k}T_{n,k+\nu}=\binom{n-k+1}{k}S_{n-k}.\hspace{3.3cm}(6)$$
 \par
 \par
 \textbf{6.} On basis of above results we can now evaluate $I_{n}^{A_{n}\oplus
 B_{n}}$. By definition $I_{n}^{A_{n}\oplus
 B_{n}}$ is the sum of $1$'s in the incidence matrix relative to
 tableau $A_{n}\oplus
 B_{n}$. Counting them by lines and because of (4) we therefore
 have
 $$I_{n}^{A_{n}\oplus
 B_{n}}=\sum_{1\leq i,j \leq S_{n}}\delta (J_{i}^{n},J_{j}^{n})
 =\sum_{i=1}^{S_{n}}\Big ( \sum_{j=1}^{S_{n}}\delta
(J_{i}^{n},J_{j}^{n})\Big )$$
$$\hspace{4cm} = \sum_{i=1}^{S_{n}}\Big(
\sum_{\nu =1}^{M(J_{i}^{n})}(-1)^{\nu
-1}\binom{M(J_{i}^{n})}{\nu}S_{n-\nu}\Big ). \hspace{2.6cm}(7)$$
Since there are $T_{n,k}^{A_{n}\oplus
 B_{n}}$ $n$-iterates with multiplicity $k$, $k=1,2,\dots,$ the
 double sum can be rearranged by pooling together all
 terms with $M(J_{i}^{n})=k$. As a consequence, using our convention, we have
 $$ \hspace{2.6cm}I^{A_{n}\oplus B_{n}}=\sum_{k=1}^{\infty}T_{n,k}^{A_{n}\oplus
B_{n}}\Bigl[
\sum_{\nu=1}^{\infty}(-1)^{\nu-1}\binom{k}{\nu}S_{n-\nu} \Bigr].
\hspace{2.6cm}(8)$$ Substituting $T_{n,k}^{A_{n}\oplus B_{n}}$ by
its value from (4) we obtain
$$I_{n}^{A_{n}\oplus B_{n}}=\sum_{k=1}^{\infty}
\frac{n+2}{k}T_{n-1,k-1}\Bigl[
\sum_{\nu=1}^{\infty}(-1)^{\nu-1}\binom{k}{\nu}S_{n-\nu}
\Bigr],$$ and reversing the order of summation we get,
$$I_{n}^{A_{n}\oplus B_{n}}=(n+2)\Big\{
S_{n-1}\Bigl[\sum_{k=1}^{\infty}\frac{1}{k}\binom{k}{1}T_{n-1,k-1}\Bigr]-
S_{n-2}\Bigl[\sum_{k=1}^{\infty}\frac{1}{k}\binom{k}{2}T_{n-1,k-1}\Bigr]+
\dots$$
$$\hspace{2.5cm}+(-1)^{n-\nu}S_{n-\nu}\Bigl[\sum_{k=1}^{\infty}\frac{1}{k}\binom{k}{\nu}T_{n-1,k-1}\Bigr]+
\dots \Bigr\}.\hspace{3cm}(9)$$
%--------------------------end page 4 ----------------------
Observing that
$\frac{1}{k}\binom{k}{\nu}=\frac{1}{\nu}\binom{k-1}{\nu-1} $ and
setting $\mu = k-1$ as a new running index (9) becomes
$$I_{n}^{A_{n}\oplus
B_{n}}=(n+2)\Bigl\{\frac{S_{n-1}}{1}\Bigl[\sum_{\mu=0}^{\infty}\binom{\mu}{0}T_{n-1,\mu}\Bigl]
-
\frac{S_{n-2}}{2}\Bigl[\sum_{\mu=0}^{\infty}\binom{\mu}{0}T_{n-1,\mu}\Bigl]+
\dots $$
$$\hspace{2.5cm}+(-1)
^{\nu-1}\frac{S_{n-\nu}}{\nu}\Bigl[\sum_{\mu=0}^{\infty}\binom{\mu}{\nu-1}T_{n-1,\mu}\Bigl]+\dots\Bigr\}.\hspace{3cm}(10)$$
The expressions in brackets can be evaluated from (6) by
replacing $n$ by $n-1$, $\binom{k+\nu}{\nu}$ by
$\binom{k+\nu}{k}$ and $k$ successively by $0, 1, \dots$. This
gives
$$k=0\hspace{1cm}
\sum_{\nu=0}^{\infty}\binom{\nu}{0}T_{n-1,\nu}=\binom{n}{0}S_{n-1}\hspace{1.8cm}
$$
$$k=1\hspace{1cm}
\sum_{\nu=0}^{\infty}\binom{\nu+1}{1}T_{n-1,\nu+1}=\binom{n-1}{1}S_{n-2}$$
$$\dots \hspace{3cm}\dots\hspace{3cm}\dots$$
$$\hspace{1.3cm}k=k\hspace{1cm}
\sum_{\nu=0}^{\infty}\binom{\nu+k}{k}T_{n-1,\nu+k}=\binom{n-k}{k}S_{n-k-1}\hspace{0.8cm}$$
$$\dots \hspace{3cm}\dots\hspace{3cm}\dots$$
 Inserting these values in (10) we obtain
$$\hspace{2cm}I_{n}^{A_{n}\oplus
B_{n}}=(n+2)\sum_{\nu=0}^{\infty}(-1)
^{\nu}\frac{1}{\nu+1}\binom{n-\nu}{\nu}S_{n-\nu-1}^{2}.\hspace{2.2cm}(11)$$
Actually this sum is finite since for $\nu>\frac{n}{2}$ all
coefficients $\binom{n-\nu}{\nu}$ are zero. It can be
asymptotically evaluated by the same heuristic  method we used in
MS  2.6 to evaluate $I_{n}$. Setting $\nu=k-1$ and after obvious
transformations, (11) becomes
$$\hspace{1.5cm}I_{n}^{A_{n}\oplus
B_{n}}=(n+2)S_{n}^{2}\,\sum_{k=1}^{\infty}(-1)
^{k-1}\frac{1}{k}\binom{n-k+1}{k-1}\Big(\frac{S_{n-k}}{S_{n}}\Big)^{2}\hspace{1.5cm}(12)$$
On the other hand for $n\rightarrow \infty$ and $k$ finite
$$\frac{S_{n-k}}{S_{n}}\backsim \frac{1}{4^{k}}$$
$$\frac{1}{k}\binom{n-k+1}{k-1}\backsim \frac{n^{k-1}}{k!}$$
so that the general term of the series in (12) behaves for
$n\rightarrow\infty$ like
$$(-1)
^{k-1}\frac{n^{k-1}}{k!}\Big(\frac{1}{4^{k}}\Big)^{2}=(-1)^{k-1}\frac{1}{k!}\frac{1}{4^{2k}}n^{k-1}.$$
%-------------------------------------end page 5---------------------------
We now sum over $k$ to get
$$I_{n}^{A_{n}\oplus
B_{n}} \backsim
(n+2)S_{n}^{2}\,\sum_{k=1}^{\infty}(-1)^{k-1}\frac{1}{4^{2}k!}\frac{1}{4^{2k-2}}n^{k-1}$$
$$\hspace{0.8cm}\backsim \frac{(n+2)S_{n}^{2}}{4^{2}}\,\sum_{k=1}^{\infty}(-1)^{k-1}\frac{1}{k!}\Big(\frac{n}{4^{2}}\Big)^{k-1}$$
$$\hspace{0.4cm}\backsim \frac{(n+2)S_{n}^{2}}{n}\,\sum_{k=1}^{\infty}(-1)^{k-1}\frac{1}{k!}\Big(\frac{n}{4^{2}}\Big)^{k}$$
$$\hspace{3.8cm}\backsim\frac{n+2}{n}\,S_{n}^{2}\big( 1-e^{-\frac{n}{16}} \big)\hspace{4.8cm}(13)$$
\par
As said in MS  2.6 above argument is not rigorous, but can be
made so by estimating for $n\rightarrow\infty$ the differences
$|\frac{S_{n-k}}{S_{n}}-\frac{1}{4^{k}}|$ and
$|\frac{1}{k}\binom{n-k+1}{k}-\frac{n^{k-1}}{k!}|$. A stronger
result can be obtained if we apply the Sterling formula to the
binomial coefficients, taking into account that even
$S_{n}=\frac{1}{n+1}\binom{2n}{n}$ contains binomial coefficients
(see SP).
\par
\textbf{6.}$\;\;I_{n}^{A_{n}\oplus B_{n}}$ was defined as the
number of formally reducible identities $J_{i}^{n}=J_{j}^{n}$ of
order $n$, so that the number of formally irreducible identities
is $S_{n}^{2}-I_{n}^{A_{n}\oplus B_{n}}$, since the total number
of $n$-identities is $S_{n}^{2}$.
\par
(13) means that the order of
$I_{n}^{A_{n}\oplus B_{n}}$ for $n\rightarrow\infty$ is
$O\big(\frac{n+2}{n}(1-e^{-\frac{n}{16}})S_{n}^{2}\big)$, which
gives for $S_{n}^{2}-I_{n}^{A_{n}\oplus B_{n}}$ the order
$O\Bigl(|\frac{n+2}{n}\;e^{-n/16}-\frac{2}{n}|S_{n}^{2}\Bigr).$}
\par
Summarizing and using the terminology of algebras (see BO), we
have proved following.
\par
\textbf{Theorem.} The number of algebras defined by a binary
operation satisfying a formally irreducible identity between two
$n$-iterates of the operation is, for $n\rightarrow\infty$,
$O\Bigl(|\frac{n+2}{n}\;e^{-n/16}-\frac{2}{n}|S_{n}^{2}\Bigr).$
\par
\textbf{Acknowledgment.} I wish to thank Peter Krikelis,
University of Athens, Dep. of Mathematics, for his help.
\par
\textbf{References.}
$$\left. \begin{matrix}\text{1. MS}\hfill &\text{Petridi, Constantin M., Krikelis, Peter
B.: Mathematical Structures}\hfill \cr &\text{ defined by
identities, arXiv: math.RA/0110333 v1 Oct 2001}\hfill  \cr
\text{2. BO}\hfill &\text{Bakhturin Yu.A., Ol'Shanskij A.Yu.: II
Identities, Encyclopedia of} \hfill  \cr & \text{Mathematical
Sciences Vol. 18, Springer - Verlag
 Berlin Heideberg 1991.} \hfill \cr
 \text{3. SP}\hfill
&\text{Stanica, Pantelimon: Good lower and upper bounds on
binomial} \hfill \cr &\text{coefficients, Journal of Inequalities
in Pure and Applied Mathematics,}\hfill \cr &\text{ Vol 2, Issue
3, Article 30, 2001}\hfill
  \end{matrix}\right. $$
\end{document}